\documentclass[a4paper, 12pt]{article}

\usepackage{amssymb,amsbsy,amsmath,amsfonts,amssymb,amscd}

\newtheorem{theorem}{Theorem}[section]
\newtheorem{lemma}[theorem]{Lemma}
\newtheorem{proposition}[theorem]{Proposition}

\newcommand{\ds}{\displaystyle}
\newcommand{\GG}{\mathbb G}
\newcommand{\ofr}{{\mathfrak o}}
\newcommand{\pfr}{{\mathfrak p}}

\newcommand{\lra}{\longrightarrow}
\newcommand{\noi}{\noindent}
\newcommand{\Ap}{\mathcal A}

\newcommand{\St}{{\mathbf S}{\mathbf t}}

\newcommand{\ZZ}{\mathbb Z}
\newcommand{\RR}{\mathbb R}
\newcommand{\CC}{\mathbb C}

\newcommand{\HH}{\mathcal H}

\newcommand{\Ch}{\rm Ch}

\title{Distinction of the Steinberg representation}
\author{P. Broussous\\
 With an appendix by F. Court\`es}

\date{\today}

\begin{document}
\maketitle

\begin{abstract}
 We prove  Dipendra Prasad's conjecture on distinction of the
 Steinberg representation [Pr] for symmetric spaces of the form $\GG (E)/\GG
 (F)$, when $\GG$ is a split reductive group defined over $F$, and
 $E/F$ an unramified quadratic extension of non-archimedean local fields.

\end{abstract}

\section*{Introduction}

 Let $\GG$ be a connected reductive group defined over a
 non-archimedean local field $F$, and let $E/F$ be a quadratic galois
 extension of $F$. If $\pi$ is a smooth representation of $\GG (E)$
 and $\chi$ a smooth character of $\GG (F)$, one says that $\pi$ is
 $\chi$-distinguished if the intertwining space
$$
{\rm Hom}_{\GG (F)}\, (\pi , \chi )
$$
is non-trivial.
\medskip

 Let $\St_E$ denote the Steinberg representation of $\GG (E)$. In
 [Pr], Dipendra Prasad defines an explicit quadratic abelian character $\chi_F$
of $\GG (F)$ and makes the following conjecture.
\bigskip

\noi {\bf Conjecture}. ([Pr], Conjecture 3, page 77). Assume that the 
derived subgroup of $\mathbb G$ is quasi-split. Then:
\medskip

(a) {\it The Steinberg
representation of $\GG (E)$ is $\chi_F$-distinguished.}
\smallskip

(b) {\it For any other smooth character $\chi$ of $\GG (F)$, different
  from $\chi_F$, the Steinberg representation of $\GG (E)$ is not
  $\chi$-distinguished.}
\bigskip

 This conjecture is proved for ${\rm GL}(n)$ (by Prasad [Pr2] when $n=2$, and
 by Anandavardhan and Rajan [AR], Theorem 1.5,  for any $n$ and without restriction
 on the quadratic field extension $E/F$).
 \medskip

 In this article, we first prove the following result.
\bigskip

\noi {\bf Theorem 1.} {\it Assume that}
\medskip

 (i) {\it $E/F$ is unramified,}
\smallskip

 (ii) {\it the residue field $k_F$ of $F$ is large enough.}
\smallskip

 (iii) the algebraic group $\GG$ is split over $F$, and to make our proof less technical:
\smallskip

 (iv) The root system of $\GG$ relative to any maximal split torus is irreducible.

\medskip

\noi {\it  Then there exists an explicit quadratic character $\epsilon_F$ of
${\mathbb G}(F)$, such that $\St_E$ is $\epsilon_F$-distinguished.}
\bigskip

 We think that conditions (ii) and (iv)  are not  necessary. On the other
 hand, conditions (i) and (iii)  are   crucial for our proof.
\medskip

 Hence,  in a particular
 case,  we obtain a proof of
 part (a) of Prasad's conjecture modulo the fact that $\epsilon_F
 =\chi_F$. This equality is true for ${\rm GL}(n)$ and when $\GG$ is
 simply connected (in this case $\chi_F =\epsilon_F =1$).
 We expect it to be always true.

\bigskip

The idea of the proof is to use the model
of the Steinberg representation as a space of harmonic functions on the
chambers of $X_E$, the building of $\GG (E)$. The building $X_F$ of
$\GG (F)$ embeds in $X_E$ as a sub-simplicial complex of same
dimension.  The $\GG (F)$-equivariant linear form is then simply
 the "period" obtained by
summing a fonction over the sub-building. The Iwahori-spherical vector is a
test vector of this linear form. The difficult point is to prove that the
restriction of a harmonic function to $X_F$ is $L^1$.
\bigskip

 We then prove the following.
\bigskip

\noi {\bf Theorem 2}. {\it Assume that assumptions (i), (iii) and (iv)
  of Theorem 1 hold. Let ${\mathbb G}^{\rm der}$ be the derived
  group of $\mathbb G$. We have the multiplicity $1$ result:}
$$
{\rm Dim}_{\mathbb C}\, {\rm Hom}_{{\mathbb G}^{\rm der}(F)}\, (\St_E
,{\mathbb C})\leq 1\ .
$$
\bigskip

 As a consequence, under the assumptions of Theorem 1, points (a) and
 (b) of Prasad's conjecture hold.
\bigskip

 Theorem 2 is a consequence of a transitivity property of the action
 of ${\mathbb G}(F)$ on the chambers of $X_E$. The proof of this
 property is provided by F. Court\`es in an appendix to this article.

\bigskip

 Since the Steinberg representation factors through $\GG (E) /Z_E$, where
 $Z_E$ is the center of $\GG (E)$, we will assume that the group $\GG$
 is semi-simple.
\bigskip

 This work was written while the first author  was supported by the french ANR grant JIVARO.
He wants to thank Fran\c cois Court\`es, Nadir Matringe and Dipendra Prasad for their help in writting the manuscript.

\section{Notation: groups and buildings}

 We fix a locally compact non-archimedean and non-discrete field $F$.
 We do not assume that the (residue) characteristic of $F$ is not
 $2$. We let $E/F$ be an unramified quadratic extension of $F$.
\medskip

 If $K$ is any locally compact non-archimedean and non-discrete field,
 we denote by
\medskip

 -- $\ofr_K$  the ring of integers of $K$,
\smallskip

 -- $\pfr_K$ the maximal ideal of $\ofr_K$,
\smallskip

 -- $k_K=\ofr_K /\pfr_K$ the residue field,
\smallskip

 -- $q_K =\vert k_K \vert$ the cardinal of $k_K$.
\medskip

 We in particular have $q_E =q_F^2$.
\bigskip

 We fix a connected semisimple  group $\mathbb G$ split and  defined
 over $F$. We  denote by $d$ its rank and by
$G=G_F$ its group of $F$-rational points. For simplicity, We shall 
 assume that the root
 system of $\mathbb G$ is irreducible.
\medskip

 We fix a maximal split torus $\mathbb T$ of $\GG$ and we  denote by
 $N$ the normalizer of ${\mathbb T}(F)$ in $G$. Let $T^0$ be the
 subgroup of $T$ generated by the $\xi (u)$, 
where $\xi$ runs over the rational cocharacters of $\mathbb T$ and
$u$  over $\ofr_F^\times$. Then 
$T^0$ is the maximal compact subgroup of $T$.
\medskip

 Let $X=X_F$ be the semi-simple Bruhat-Tits building of $\GG$. 
This is a locally compact topological space on which $G$ acts
continuously. It has dimension $d$.  The space $X$ is naturally the 
geometric realization of a  simplicial complex and the group $G$ 
acts by preserving the simplicial structure.

  Let us fix a chamber $C_0$ in the apartment $A$ of $X$ attached 
to $T$ and write $I$ for the Iwahori subgroup of $G$ attached to
$C_0$.  This is the pointwise stabilizer of $C$ in $G$.  

 By [IM] (also see [I]), the affine Weyl group $W=N(T)/T^0$ of $T$ 
may be written as a semidirect 
product $W=\Omega\ltimes W_0$ of a coxeter group 
$W_0$ by a finite abelian group $\Omega$,   in such a way that:

\medskip

 (a) $\Omega$ normalizes $I$,
\smallskip

 (b) if $N^0$ is the inverse image of $W_0$ in $N$,
 $(I,N_0 )$ is a Tits system (or $BN$-pair),
\smallskip

 (c) the set
$$
G_0 =IW_0 I =\bigcup_{w\in W_0}IwI
$$
\noi is a normal subgroup of $G$ and $G/G_0 \simeq \Omega$.
\medskip

 \noi The pair $(I,N)$ is a {\it generalized Tits system}. When $\GG$ is
 simply connected, we have $\Omega =\{ 1\}$. 
\medskip

 As a simplicial complex $X_F$ is the building of the $BN$-pair
 $(I,N_0 )$ [BT].  In particular  $X_F$ is labellable (in the sense of
 [Br], Appendix C, page 29) and $G_0$ acts
 on  $X_F$ by preserving the labelling of simplices.

\bigskip

The group $\Omega$ acts on $A_0$ and stabilizes the chamber $C_0$.
 For $\omega\in \Omega$, we  denote by $\epsilon (\omega )$ the
 signature of the permutation induced by the action of 
$\omega$ on the vertex set of $C_0$. We define a quadratic character
$\epsilon_{G_F}$ of $G_F$ by 
$$
\epsilon_{G_F} =\epsilon \circ p_0
$$
\noi where $p_0$~: $G\lra G/G_0$ denotes the canonical projection.
\bigskip

We fix an unramified quadratic extension $E/F$. By [T] there 
is a canonical embedding 
$$
j~:\ X_F\lra X_E
$$
\noi  of $X_F$ in the semisimple building $X_E$ of $\mathbb G$ over $E$. The
 Galois group ${\rm Gal}(E/F)$ acts on $X$ and  $j$ 
is ${\rm Gal}(E/F)\ltimes G_E$-equivariant. Moreover since $\mathbb G$
is split and $E/F$ is 
unramified, we have that :
\medskip

 -- $j(X_F )$ is the set of ${\rm Gal}(E/F)$-fixed points in $X_E$,
\smallskip

 -- $j$ is simplicial.
\smallskip

 --  $X_F$ and $X_E$ share the same 
dimension $d$, and $j$ maps  chambers to chambers. 
\medskip

 We shall identify $X_F$ as a subsimplicial complex of $X_E$
by viewing $j$ as an inclusion. 
\bigskip

 Let $\Delta_d$ be the standard abstract simplex of dimension $d$.
 We view its set of simplices as the power set of 
$\lbrace 0,1,...,d\rbrace$.  Since $X_E$ is labellable, there exists a
simplicial map
$$
\lambda_E~: \ X_E\lra \Delta_d
$$
\noi which preserves the dimension of simplices. If $\sigma$ is
a simplex of $X_E$, we call  $\lambda_E (\sigma )$ its {\it label}
or {\it type}. The restriction $\lambda_F =(\lambda_E )_{\vert X_F}$ is a labelling
of $X_F$ preserved by the action of $G_0$.

  Let $C$ be a chamber of $C$ and  $g\in G_E$.
  Let $(s_0 ,...,s_d )$ (resp. $(t_0 ,...,t_d )$) be an ordering
of the vertices of  $C$ (resp. of $gC$) such that
$\lambda_E (s_i )=\{ i\}$, $i=0,...,d$ (resp. $\lambda_E (t_i )=i$,
$i=0,...,d$). We denote by $\epsilon (g,C)$ the signature of the 
permutation:
$$
\left(
\begin{array}{cccc}
g.s_0 & g.s_1 & \dots & g.s_d \\
t_0   & t_1   & \dots & t_d
\end{array}
\right)
$$

\begin{lemma} With the previous notation, we have:

(i) the signature
$\epsilon (g,C)$ does not depend on $C$.

(ii) The map
$g\mapsto \epsilon (g)=\epsilon (g,C_0 )$ is a character
of $G_E$.

(iv) The character $\epsilon$  satisfies $\epsilon_{\vert G_F}= \epsilon_{G_F}$.
\end{lemma}

\noi {\it Proof}. It is easy and based on the fact that the subgroup $G_E^{*}$ of $G_E$, 
formed of those elements preserving the labelling $\lambda_E$, acts transitively on chambers of $X_E$.
 Details are left to the reader.  $\square$

\section{The Steinberg representation}

 There are several equivalent definitions of the Steinberg 
representation $\St_E$ of $G_E$. That we shall use comes
from the following beautiful theorem due to Borel and Serre.

\begin{theorem} [BS] The representation of $G_E$ in $H_c^d  (X_E ,\CC  )$, 
the $d$-th cohomology space with compact support, 
with coefficient in $\CC$, where $d$ is the $E$-rank of $\GG$,
 is equivalent to the Steinberg representation.
\end{theorem}

 Let $\Ch_E$ denote the set of chambers of $X_E$ and $\CC [\Ch_E ]$
the $\CC$-vector space of complex valued fontions on $\Ch_E$ of
arbitrary support. A function $f\in \CC [\Ch_E ]$ is
called a {\it harmonic cocycle} if for all codimension $1$ 
simplex $D$ of $X_E$, we have
$$
\sum_{C\supset D} f(C)=0
$$
\noi where the sum is over the chambers of $X_E$ that contain $D$ as 
a subsimplex. We denote
by $\HH (X_E )$ the $\CC$-vector space of harmonic cocycles on $X_E$. 

 We define a linear representation $(\pi_E ,\HH (X_E ))$ of $G_E$ 
in $\HH (X_E )$ by the formula:
$$
[\pi_E (g).f](C)=\epsilon (g) f(g^{-1}C)\ , \ g\in G ,\ C\in \Ch_E\ .
$$
\noi This representation is not smooth in general and we denote by
 $(\pi_E ,\HH (X_E )^{\infty})$ its smooth part. 

\begin{proposition} The representation $(\pi_E ,\HH (X_E )^{\infty})$
is equivalent as a $G_E$-representation to the contragredient of $\St_E$.
\end{proposition}
 
\noi {\it Proof}. For $k=d-1,d$, let $C_c^k (X_E )^{\rm alt}$ be the
$\CC$-vector space of alterned $k$-cochains on $X_E$ with coefficients
in $\CC$, the field of complex numbers. Denote by $\Ch_E^*$ the set of
pairs $(C,\sigma)$ formed of a chamber $C$ of $X_E$ together with a
bijection $\sigma$ from the vertex set of $C$ to $\{
0,1,...,d\}$. We let $G_E$ act on $\Ch_E^*$ by $g.(C,\sigma
)=(g.C,\sigma \circ g^* )$, where $g^*$ is the bijection from the
vertex set of $C$ to the vertex set of $g.C$ induced by $g$. Then
$C_c^d (X_E )^{\rm alt}$ is the set of maps $f$~: $\Ch_E^*\lra \CC$
satisfying:
\medskip

 -- $f$ has finite support,
\smallskip

 -- for all $(C, \sigma )\in \Ch_E^*$ and for all permutation $\tau$
 of $\{ 0,...,d-1\}$, we have
$$
f(C, \tau\circ \sigma )= \epsilon (\tau )f(C,\sigma )
$$
\noi where $\epsilon (\sigma )$ denotes the signature of $\sigma$.
\medskip

 The group $G_E$ naturally acts on $C_c^d (X_E )^{\rm alt}$. Similarly
 we define the $G_E$-module $C_c^{d-1}(X_E )^{\rm alt}$. The
 coboundary map 
$$
d~: \ C_c^{d-1}(X_E )^{\rm alt}\lra  C_c^{d}(X_E )^{\rm alt}
$$
\noi is given by 
$$
dh \, (C ,\sigma )=\sum_{D\subset C} h(D, \sigma_{\vert D})\ ,
\ (C,\sigma )\in \Ch_E^*
$$
\noi where $\sigma_{\vert D}$ denotes the restriction of $\sigma$ to
the vertex set of $D$. 

 \medskip

 For $k=d-1,d$, let $C_c^k (X_E )$ be the $\CC$-vector space of usual
 $k$-cochains with finite support. By orienting the simplices of $X_E$
 thanks to the labelling $\lambda$, we obtain a coboundary map $d$~:
 $C_c^{d-1}(X_E )\lra C_c^{d}(X_E)$ given by
$$
dh (C)=\sum_{D\subset C} (-1)^{\lambda (C\setminus D)} h(D) .
$$ 
For $k=d-1,d$, we have an isomorphism of $\CC$-vector spaces
$$
C_c^{k} (X_E )^{\rm alt} \lra C_c^{k}(X_E )
$$
\noi given by 
$$
f \mapsto \lbrace C \mapsto f(C,\lambda_{\vert C})\rbrace
$$
\noi where $\lambda_{\vert C}$ denotes the restriction of the
labelling $\lambda$ to the vertex set of the simplex $C$. These
isomorphisms are $G_E$-equivariant if one lets $G_E$ act on $C_c^k
(X_E )$ via
$$
[g.f](D) = \epsilon_{G_E} (g)\,  f(g^{-1}.D) ,\ D \ k\text{-simplex of
} X_E \ .
$$
\noi Moreover the isomorphisms are compatible with the coboundary
maps.
\medskip

 The space $H_c^d (X_E )$ is known to be isomorphic to $C_c^{d}(X_E
 )^{\rm alt} /dC_c^{d-1}(X_E )^{\rm alt}$ as a $G_E$-module. So it is
   isomorphic to $C_c^d (X_E )/dC_c^{d-1}(X_E )$ as a $G_E$-module. 

By  letting $V^*$ denote the algebraic dual of a $\CC$-vector space
$V$, we have
$$
\big( H_c^d (X_E )\big)^* =\lbrace \omega \in C_c^d (X_E )^*\ ; \ 
f_{\vert dC_c^{d-1}(X_E )}=0\rbrace\ .
$$
We may identify $ C_c^d (X_E )^*$ with $\CC [\Ch_E ]$ by using the
pairing
$$
\langle \omega ,f\rangle =\sum_{C\in\Ch_E} \omega (C) \, f(C)\ ,
\ \omega \in \CC [\Ch_E ] , \ f\in C_c^{d}(X_E )\ .
$$
Then for $\omega \in \CC [\Ch_E ]$, the condition $\omega_{\vert dC_c^{d-1}(X_E
  )}=0$ writes $ \langle \omega , dh\rangle =0$, for all $h\in
C_c^{d-1} (X_E )$. This may be rewritten
$$
\langle d^* \omega ,h\rangle = 0\ , h\in C_c^{d-1}(X_E)\ \text{that is
} d^* \omega =0
$$
\noi where $d^*$~: $C^d (X_E)^* \lra C_c^{d-1}(X_E )^*$ is the adjoint
of $d$. But a simple computation shows that
$$
d^* \omega (D) = \sum_{C\supset D} \omega (C)\ , \ D\ (d-1)\text{-simplex}
$$
\noi so that $d^*\omega =0$ is the harmonicity condition. $\square$

\bigskip

 Note that the Steinberg representation of $G_E$ is self-dual.

\section{Some geometric lemmas}

 We denote by $d_g$ the combinatorial distance on $X_E$ defined as follows.
 For $C$, $D\in \Ch_E$, $d_g (C,D)$ is the length $k$ of a minimal gallery
$(D_0 , D_1 ,...,D_k )$ satisfying $D_0 =C$ and $D_k =D$. The following result, 
due to F. Bruhat, will be very useful.

\begin{lemma} (Lemma 4.1 of [Bo]) Let $U$ be a compact open subgroup of $G_E$.
 There exists an integers $k_0 = k_0 (U)$ satisfying the following property.
 For all chamber $C$ such that $d_g (C_0 ,C)\geqslant k_0$,
there exists a chamber $D$ adjacent to $C$ such that:
\medskip

 (i) $d_g (C_0 ,D)=d(C_0 ,C)-1$;
\smallskip

 (ii) the group $U$ acts transitively on the set of chambers $C'$ such that
 $C'\not= D$ and $C'\cap D =C\cap D$.
\end{lemma}

\begin{lemma} Let $D$ be a codimension $1$ simplex in $X_F$ (resp. in
  $X_E$). Then $D$ is contains in $q_F +1$ chambers of $X_F$ (resp. in
  $q_E +1$ chambers of $X_E$).
\end{lemma}

\noi {\it Proof}. We give a proof for $X_F$. Let $P_D$ be the
parahoric subgroup of $G_F$ attached to $D$ and $P_D^1$ its
pro-unipotent radical. Then $P_D /P_D^1 =\GG_D (k_F )$, where $\GG_D$
is a reductive group defined over $k_F$ and of $k_F$-rank $1$. The
chambers $C$ of $X_F$ containing $D$ are in bijection with the Borel
subgroup of $P_D /P_D^1$ by
$$
C \mapsto P_C \ \text{mod} \ P_D^1
$$
\noi where $P_C$ denotes the Iwahori subgroup attached to $C$. But
$\GG_D$ being of $k_F$-rank $1$, $\GG_D (k_F )$ posseses $q_F +1$
Borel subgroups.  $\square$
\bigskip

For any non negative integer $k$, we denote by  $\Sigma_F
(k)$ the set  of chambers of  $X_F$ at distance
$k$ from $C_0$ and set $N_k (k) =  \vert\, \Sigma_F (k)\, \vert$. 

\begin{lemma} We have
$$
N_F (k)\leqslant (d+1)\,  d^{k-1} \, q_F^{k}\ , \ k\geqslant 1\ .
$$
\end{lemma}

\noi {\it Proof}. Any chamber of $X_F$ has $d+1$ codimension $1$
faces. A chamber in $\Sigma_F (1)$ contains one of the $d+1$
codimension $1$ faces of $C_0$. By Lemma (3.2), such a face is
contained in $q_F$ chambers different from $C_0$, so that
$$
N_F (1) = (d+1)\, q_F\ .
$$
\noi Moreover, for $k\geqslant 1$, any chamber in $\Sigma_F (k)$ is
adjacent to at most $dq_K$ chambers at distance from $C_0$ greater
than $k$. The formula follows by induction on $k$. $\square$

\begin{lemma} Let $f\in \HH (X_E )^{\infty}$. There exist an integer 
$k_f$ and a positive real number $K_f$ such that the following holds.
For all $C\in \Ch_E$ such that $d_g (C_0 ,C)\geqslant k_f$, we have
$$
\vert f(C)\vert \leqslant K_f \, . \, q_E^{-d_g (C_0 ,C)} \ .
$$
\end{lemma}

\noi {\it Proof}. Since $f$ is smooth under the action of $G$, it is
fixed by an open compact subgroup $U$ small enough. Set $k_f =k_0
(U)$. For $k\geqslant 0$, set 
$$
M_k ={\rm Max}\, \lbrace \vert\,  f(C)\, \vert \ ;  \ C\in \Sigma_E
(k)\rbrace \ .
$$
\noi We are going to prove that for $k\geqslant k_f$ we have
$M_{k+1}\leqslant q_E^{-1} M_k$; the lemma will follow. 

 Let $C\in \Sigma_E (k+1)$. By applying Lemma (3.1), there exists
 $D\in \Sigma_E (k)$ such that $U$ acts transitively on 
$$
[C,D]:=\lbrace G\in \Ch_E\ ; \ G\not= D\text{ and } G\cap D = C\cap
D\rbrace\ .
$$
\noi It follows that $f$ is constant on $[C,D]$. By applying the
harmonicity condition at the codimension $1$ face $C\cap D$, we get 
$$
q_E \, f(C) +f(D) = 0\ ,
$$
since $[C,D]$ has $q_E$ elements. So $\vert \, f(C)\, \vert =q_E^{-1}
\vert\, f(D)\, \vert$, and our assertion follows. $\square$

\begin{lemma} Assume that $q_F >d$. Let $f\in \HH (X_E)^\infty$. Then
  we have
$$
f_{\vert \Ch_F} \in {\rm L}^1 \, (\Ch_F )\ 
$$
\noi where $L^1 (\Ch_F )$ denotes the set of complex functions $g$ on
$\Ch_F$ such that
$$
\sum_{C\in \Ch_F}\vert\, g(C)\, \vert < +\infty\ .
$$

\end{lemma}

\noi {\it Proof}.  We may write 
$$
\sum_{C\in \Ch_F}\vert\, f(C)\, \vert =\sum_{k\geqslant 0}\sum_{C\in
  \Sigma_F (k)} \vert\, f(C)\, \vert\ .
$$
\noi By the previous lemmas, for $k$ large enough  and for some
constant $K>0$, we have:
$$
\sum_{C\in
  \Sigma_F (k)} \vert\, f(C)\, \vert \leqslant K\,
(\frac{dq_F}{q_E})^k\ .
$$
\noi with $q_E = q_F^2$. The result follows. $\square$
 \bigskip

\noi {\it Remark}. If $\GG$ is of rank $1$, then the condition $q_F
>d$ is automatically satisfied.

\section{Constructing $G_F$-equivariant linear forms}

 In this section, we assume, as in Lemma (3.5), that we have $q_F
 >d$. 
\bigskip

 Thanks to lemma (3.5), the linear map $\lambda$ on $\HH (G_E
 )^\infty$ given by
$$
\lambda (f)=\sum_{C\in \Ch_F} f(C)
$$
\noi is well defined. For $g\in G_F$ and $f\in \HH (G_E )^\infty$, we
have
$$
\lambda (\pi_E (g).f)=\sum_{C\in \Ch_F} \epsilon_{G_F} (g)
f(g^{-1}C)=\epsilon_{G_F} (g) \lambda (f)
\ .
$$
\noi Hence we have $\lambda \in {\rm Hom}_{G_E}(\St_E ,\epsilon_{G_F}
)$.

\begin{theorem}  The Steinberg representation of $G_E$ is
  $\epsilon_{G_F}$-distinguished.

\noi  More precisely,  a non-zero 
 Iwahori-spherical vector  is a test vector for $\lambda$.
\end{theorem}

\noi {\it Proof}. It suffices to prove that $\lambda$ is not
trivial. Let $f$ be the Iwahori-spherical vector in $\HH (G_E
)^\infty$ normalized in such a way that $f(C_0 )=1$. 
 In Lemma (3.1), if $U=I$, we may take $k_0 = 0$. It follows from the
 proof of Lemma (3.4) that, for all $k\geqslant 0$, $f$ has constant
 value $\ds (\frac{-1}{q_E})^{-k}$ on $\Sigma_F (k)$. As a consequence
$$
\lambda (f)= \sum_{k\geqslant 0} \big( \sum_{C\in \Sigma_F (k)} f(C)\big)
$$
is an alternating series. In particular we have
$$
\sum_{C\in \Sigma_F (0)} f(C) > \lambda (f) >
\sum_{C\in \Sigma_F (0)} f(C)+\sum_{C\in \Sigma_F (1)} f(C)
$$
\noi that is

$$
1 > \lambda (f) > 1-\frac{d+1}{q_F}\geqslant 0\ .
$$
\noi and our Theorem follows. $\square$
\medskip

 Note  that if  the $F$-rank of $\GG$ is $1$ the value $\lambda (f)$
 may be explicitely computed. Indeed in that case,  $X_F$ is a regular
 tree of valence $q_F +1$, and we have $N_F (k) = 2q_F^k$, $k\geqslant
 1$. Hence
$$
\lambda (f)=1 +\sum_{k\geqslant 1} 2q_F^k (\frac{-1}{q_E})^k =
 1-\frac{2}{q_F +1} \ .
$$

\section{Multiplicity $1$}

In this section we release the condition $q_F >d$ and prove Theorem 2
without restriction on the size of $k_F$. 
\bigskip

 Set ${\mathbb H}={\mathbb G}^{\rm der}$ and $H={\mathbb H}(F)$. Note
 that $\mathbb H$ and $\mathbb G$ share the same (semisimple)
 Bruhat-Tits building over $F$ (resp. over $E$). This essentially
 comes from the fact that the inclusion $H\lra G$ is $B-N$-{\it
   adapt\'e} in the sense of [BT] (1.2.13), page 18 (cf. [BT]
 {\S}2.7., page  49).  By Proposition (2.2), we have 
$$
{\rm Hom}_H\, (\St_E ,\CC )={\mathcal H}(X_E )^H\ ,
$$
\noi the space of harmonic cochains on $X_E$ fixed by $H$.

\bigskip

 Our proof relies on the following fundamental result whose proof is given in the
appendix.

\begin{theorem} Let $C$ be a chamber of $X_E$ at combinatorial  distance $\delta
  \geqslant 0$ from $X_F$. Then $G_F$ acts transitively on the set
  ${\rm Ch}\, [C,\delta +1 ]$ of chambers $D$ of $X_E$ satisfying:
\medskip

 -- $D$ and $C$ are adjacent,
\smallskip

 -- $d(D,X_F )=\delta +1$.
\end{theorem}
\bigskip

 Note that since $H$ is contained in $G_0$ ($G/G_0$ is abelian), it
 acts on $\HH (X_E )$ via the formula
$$
h.\omega (C) = \omega (h^{-1}C ) \ , \ h\in H, \ \omega \in \HH (X_E
)\ .
$$

\bigskip
 Let $\omega\in {\mathcal H}(X_E )^H$. Since $H$ acts transitively on
 $X_F$, the value $\omega (C)$ does not depend on the chamber $C$ of
 $X_F$. Let us denote it by $\varphi (\omega )$. Theorem 2 is a
 consequence of the following:

\begin{lemma} The linear map
$$
\varphi~: \ {\mathcal H}(X_E )^H \lra \CC\ , \ \omega \mapsto \varphi (C)
$$ 
\noi is injective.
\end{lemma}

\noi {\it Proof}. For all integers $\delta \geqslant 0$, let  ${\rm
  Ch}\, (X_F ,\delta )$ denote the set of chambers in $X_E$ at
combinatorial distance $\delta$ from $X_F$ (in particular ${\rm Ch}\,
(X_F ,0)$ is the set of chambers of $X_F$). Let $\omega \in \HH (X_E
)^H$ and $\delta \geqslant 0$ be an integer.
 We prove that the restriction $\omega_{\vert  {\rm
  Ch}\, (X_F ,\delta +1)}$ is entirely determined by the restriction 
$\omega_{\vert  {\rm Ch}\, (X_F ,\delta)}$. The lemma will obviously
follow. 
\bigskip

 Let $D\in {\rm Ch}\, (X_F ,\delta +1)$. Fix a chamber $C\in  {\rm
   Ch}\, (X_F ,\delta )$ adjacent to $D$ and set $M=C\cap D$. The
 harmonicity condition at the codimension $1$ face $M$ writes
$$
\sum_{\Delta \in C_M} \omega (\Delta )=0\ ,
$$
\noi where $C_M$ is the set of chambers of $X_E$ containing $M$. We
may split the set $C_M$ into two subsets : $C_M^{\delta +1}:={\rm Ch}[C,\delta +1 ]$
and its complement $C_M^{\delta}$, contained in ${\rm Ch}\, (X_F ,\delta
)$. By theorem (5.1) and the $H$-invariance of $\omega$, we have
$$
\sum_{\Delta \in C_M^{\delta +1}}\omega (\Delta ) =\vert C_M^{\delta
  +1}\vert \times \omega (C)\ .
$$

\noi Hence the harmonicity condition gives
$$
\omega (C) = -\frac{1}{ \vert C_M^{\delta
  +1}\vert} \times \sum_{\Delta \in C_M^\delta}\omega (\Delta )
$$
\noi This proves that the value $\omega (C)$ depends only on the
restriction $\omega_{\vert {\rm Ch}(X_F ,\delta )}$, and we are done. $\square$.

\appendix
\section{A transitivity result}

For every facet $A\subset X_E$, we shall denote by $K_{A,E}$ (resp $K_{A,F}$)
 the connected fixator of $A$ in $G_E$ (resp. the intersection with $G_F$ of
 that connected fixator). More generally, for every subset $S$ of $X_E$,
 we shall denote by $K_{S,E}$ (resp. $K_{S,F}$) the intersection of the $K_{A,E}$ 
(resp.$K_{A,F}$), where $A$ runs over the set of facets of $X_E$ whose intersection
 with $S$ is nonempty. Let $\overline{S}$ be the closure of $S$; we have of
 course $K_{\overline{S},E}=K_{S,E}$ and $K_{\overline{S},F}=K_{S,F}$.

\begin{proposition}
Let $d$ be a nonnegative integer, and let $C$ be any chamber of $X_E$ such
 that the combinatorial distance between $C$ and $X_F$ is $d$. Let $C'$ be a
 chamber of $X_E$ neighbouring $C$ and whose combinatorial distance from $X_F$ 
is $d+1$, let $A$ be the unique facet of codimension $1$ of $X_E$ contained in
 both $\overline{C}$ and $\overline{C'}$, and let $\Delta$ be the set of chambers
 of $X_E$ containing $A$ in their closure and whose combinatorial distance from
 $X_F$ is $d+1$. Then the group $K_{C,F}$ acts transitively on $\Delta$.
\end{proposition}

\noi {\it Proof}. Assume first $d=0$, that is $C$ is contained in $X_F$.
 Let $\Ap$ be an apartment of $X_F$ containing $C$, let $T$ be the corresponding 
$F$-split maximal torus of $G_E$, let $\Phi$ be the root system of $G_E$ relatively
 to $T_E$ and let $\pm\alpha$ be the elements of $\Phi$ corresponding to the hyperplane
 $H$ of $\Ap$ containing $A$. For every $\beta\in\Phi$, let $U_\beta=U_{\beta,E}$ be the 
corresponding root subgroup of $G_E$, and let $v$ be a normalized valuation 
(that is a valuation such that for every $\beta$, $v(U_\beta)=\ZZ\cup\{\infty\}$; 
such a valuation exists because $\GG$ is split over an unramified extension of $E$)
 on the root datum $(G,T,(U_\beta))$ such that, with the subgroups $U_{\beta,i}$ of 
$U_\beta$ being defined according to that valuation, we have
 $U_{\pm \alpha}\cap K_A=U_{\pm \alpha,0}$; we shall assume $\alpha$ is the one such that
 $U_\alpha\cap K_C=U_{\alpha,1}$.

Let $\phi$ (resp. $\phi'$) be a $F$-isomorphism between $E$ and $U_{\alpha,E}$
 (resp. $U_{-\alpha,E}$) preserving the valuation, that is such that for every
 integer $i$, $U_{\alpha,i}$ (resp. $U_{-\alpha,i}$) is the image by $\phi$ 
(resp. $\phi'$) of the elements of $E$ of valuation $\geq i$; the elements
 of $\Delta$ are the chambers of the form $\phi(x)C$, where $x$ is an element 
of $\ofr_E$ which belongs neither to $F$ nor to $\pfr_E$; moreover, $\phi(x)C$ 
depends only of the class of $x$ modulo $\pfr_E$. we can thus label the elements
 of $\Delta$ as $C_x=\phi(x)C$, where $x$ is an element of $k_E-k_F$.

Let now $\Phi^\vee$ be the system of coroots of $T_E$ associated to $\Phi$, 
and let $\alpha^\vee$ be the $1$-parameter subgroup of $T_E$ corresponding
 to $\alpha$ in $\Phi^\vee$; for every $y\in\ofr_E^*$, we have $\alpha^\vee(y)C=C$,
 and if $y$ is an element of $\ofr_F^*+\pfr_E$, $\alpha^\vee(y)$ permutes
 the elements of $\Delta$. Moreover, $\alpha^\vee(y)$ depends only of the class
 of $y$ modulo $1+\pfr_E$, hence we can view $y$ as an element of $k_E^*$.

Let $x$ be an element of $k_E-k_F$ (arbitrarily fixed for the moment). 
For every $a\in k_F^*$ and every $b\in k_F$, we have:
\[\phi(b)\alpha^\vee(a))C_x=\phi(b)(Ad(\alpha^\vee(a))\phi(x))C=\phi(b)\phi(a^2x)C=C_{a^2x+b}.\]
Hence for every element of $k_E-k_F$ of the form $y=a^2x+b$, $C_y$ is in the $G_F$-orbit of $C_x$.
 If $char(k_E)=2$, every element of $k_F$ is a square, and since $(1,x)$ is a basis of
 the $k_F$-vector space $k_E$, every element of $k_E-k_F$ is of that form, which proves 
the proposition in that case.

Now assume $p\neq 2$; there exists then $x\in k_E-k_F$ such that $\frac 1c=x^2\in k_F$; $c$ 
is then not a square in $k_F$. Let $D$ be the chamber of $\Ap$ such that
 $\overline{D}\cap\overline{C}=\overline{A}$; we have $D=nC$, where $n$ is any 
representative in the normalizer of $K_T$ in $K_A$ of the element $s_\alpha$ in the 
Weyl group of $\Phi$. Moreover, according to \cite[6.1.3 a) and b)]{BT}, we can 
assume that every such element is of the form:
\[n=\phi'(y)\phi(-y^{-1})\phi'(y),\]
with $y\in\ofr_E^*$. We then have:
\[\phi'(y)D=\phi'(y)\phi'(-y)\phi(y^{-1})\phi'(-y)C=C_{y^{-1}},\]
since $\phi'(y)C=C$. Hence $C_x=\phi(x^{-1})D=\phi(xc)D$. By the same reasoning 
as above, for every $a\in k_F^*$ and every $b\in K_F$, $\phi'(a^2xc+b)D=C_{\frac 1{a^2xc+b}}$. On
 the other hand, we have:
\[\frac 1{a^2xc+b}=\frac{a^2xc-b}{(a^2xc+b)(a^2xc-b)}=\frac{a^2xc-b}{a^4c-b^2}=\frac{x-\frac b{a^2c}}{a^2-\frac{b^2}{a^2c}}.\]

On the other hand, it is well-known and easy to check that there exists $a,b$ such that $a^2-\frac{b^2}{a^2c}$ 
is not a square; we thus obtain that there exists $a',b'$, such that $a'$ is not 
a square and $C_{a'x+b'}$ is in the $K_{C,F}$-orbit of $C_x$. By the same reasoning as 
above once again, we obtain that it is true for every $C_{a'a^2x+b+b'}$, $a\in k_F^*$, 
$b\in k_F$. Since $(k_F^*)^2$ is of index $2$ in $k_F$, we finally obtain that all of 
the $C_x$, $x\in k_E-k_F$, are in the same $K_{C,F}$-orbit, which completes the proof
 of the proposition when $C\subset X_F$.

Now assume $d>0$. Set $\Gamma=\rm{Gal}(E/F)$, and let $\gamma$ be the unique nontrivial 
element of $\Gamma$. First we prove the following lemma:

\begin{lemma}
There exists a $\Gamma$-stable apartment of $X_E$ containing both $C$ and $\gamma(C)$.
\end{lemma}

\noi {\it Proof}. Let $\Ap$ be any apartment of $X_E$ containing both $C$ and $\gamma(C)$;
 such an apartment exists by \cite[proposition 2.3.1]{BT}. Obviously, $\gamma(\Ap$)
 satisfies the same property; there exists then $g\in G_E$ such that $g\Ap=\gamma(\Ap)$,
 and we can assume $g\in K_{C,E}\cap K_{\gamma(C),E}$.  The element $\gamma(g)$ then also
 belongs to $K_{C,E}\cap K_{\gamma(C),E}$, and we have $\gamma(g)\gamma(\Ap)=\Ap$. Hence 
$\gamma(g)g$ fixes $\Ap$ pointwise, which means that it belongs to the unique parahoric
 subgroup $K_T$ of the $E$-split maximal torus $T$ of $G_E$ associated to $\Ap$. 

Let now $F_{nr}$ be the maximal unramified extension of $F$, let $G_{F_{nr}}$ be the group
 of $F_{nr}$-points of $\GG$, and let $K_{C,F_{nr}}$ be the connected fixator of $C$ viewed
 as a chamber of the Bruhat-Tits building $X_{F_{nr}}$ of $G_{F_{nr}}$. By \cite[lemma 5.1]{Cou}, 
there exists an element $h\in K_{F_{nr}}$ such that $g=\mathbf{F}(h)^{-1}h$, with $\mathbf{F}$
 being the Frobenius element of $\rm{Gal}(F_{nr}/F)$. Moreover, the restriction of $\mathbf{F}$
 to $E$ is $\gamma$, and we have:
\[\gamma(g)g=\mathbf{F}^2(h)^{-1}h\in K_T,\]
Let $T_{nr}$ be the maximal torus of $G_{F_{nr}}$ associated to $\Ap$, and let $K_{T_{nr}}$ be
 its unique parahoric subgroup; we have $K_T=K_{T_{nr}}\cap G_E$. Moreover, the Frobenius
 element of $\rm{Gal}(F_{nr}/E)$ is $\mathbf{F}^2$; by \cite[lemma 5.1]{Cou} again, there
 exists then $t\in K_{T_{nr}}$ such that $\gamma(g)g=\mathbf{F}^2(t)t^{-1}$.  Hence
 $ht=\mathbf{F}^2(ht)$, which simply means that $ht\in G_E$. We finally obtain:
\[ht\Ap=h\Ap=\gamma(h)\gamma(\Ap)=\gamma(h\Ap)=\gamma(ht\Ap),\]
hence $ht\Ap$ is a $\Gamma$-stable apartment of $X_E$ containing both $C$ and $\gamma(C)$
 and the lemma is proved. $\Box$.

Now we designate by $\Ap$ the apartment given by the above lemma, and by $T$ the
 corresponding $E$-split maximal torus of $G_E$; $T$ is defined over $F$, but not 
$F$-split. Let $\Phi$ be the root system of $G_E$ relatively to $T$, and let
 $\alpha\in\Phi$ be defined as in the case $d=0$. Since $T$ is defined over 
$F$, $\Gamma$ acts on $\Phi$.

Let $D$ be the unique chamber of $\mathbf{A}$ such that
 $\overline{C}\cap\overline{D}=\overline{A}$. Since $\Delta$ is nonempty, 
the combinatorial distance between $D$ and $\mathbf{B}_F$ must be either $d$ or $d+1$.

\begin{lemma}
Assume $H=\gamma(H)$. Then the combinatorial distance between $D$ and $\mathcal{F}$ is $d$.
\end{lemma}

\noi {\it Proof}. Let $s_H$ be the orthogonal reflection on $\Ap$ whose kernel is $H$. 
Since $H=\gamma(H)$, $\gamma$ and $s_H$ commute, hence there exists $g_H\in G_F$ such 
that $g_H$ acts on $\Ap$ via $s_H$. Let $C=C_0,\dots,C_d$ be a minimal gallery of length 
$d$ between $C$ and some chamber $C_d$ of $X_F$. Then $D=g_HC_0,\dots,g_HC_d$ is also a 
minimal gallery and $g_HC_d\subset X_F$, hence the combinatorial distance between $D$
 and $X_F$ is at most $d$. The other inequality follows from the above remarks. $\Box$

Note that the fact that $H=\gamma(H)$ implies in particular that $\gamma(\alpha)=\pm\alpha$.
 Conversely, we have:

\begin{lemma}
Assume $\gamma(\alpha)=\alpha$. Then $H=\gamma(H)$.
\end{lemma}

\noi {\it Proof}. Let $\overline{\alpha}$ be the affine root of $T$ corresponding
 to $H$; it is an affine linear form on the affine space $\Ap$, and the corresponding
 linear form on the vector space $(X_*(T)/X_*(Z))\times\RR$, where $Z$ is the center
 of $G$, is $\alpha$. Hence $\gamma(\overline{\alpha})$ is of the form
 $\overline{\alpha}+c$, with $c$ being some constant. We then have
 $\gamma^2(\overline{\alpha})=\overline{\alpha}+2c$; since $\gamma^2$ is 
trivial, it implies $c=0$, hance $H=\gamma(H)$. $\Box$

Note that it is not true when $\gamma(\alpha)=-\alpha$.

Now we prove the proposition when $H=\gamma(H)$. Consider the rank $1$ subgroup
 $G_\alpha$ of $G_E$ generated by $T$, $U_\alpha$ and $U_{-\alpha}$; it is defined over
 $F$, and the fact that $H=\gamma(H)$ implies that $G_\alpha\cap K_A=G_\alpha\cap K_{\gamma(A)}$,
 hence $K_A$ is $\Gamma$-stable. The elements of $\Delta$ are then of the form $uC$, where 
$u$ is an element of $U_\alpha$ not belonging to $G_F$, and we can finish the proof the same
 way as in the case $d=0$.

Assume now $H\neq\gamma(H)$. Let $\mathcal{C}$ be the connected component of 
$\Ap-(H\cup\gamma(H))$ containing $C$. Assume $\mathcal{C}$ contains $\gamma(C)$
 as well. Consider an apartment of $X_E$ of the form $\phi(x)\Ap$, where $\phi$ 
is defined as in the case $d=0$ for a given normalized valuation $v$ on $(G,T,(U_\beta))$, 
and $x$ is an element of $E$ of valuation $i$, where $i$ is such that 
$U_\alpha\cap K_{C,E}=U_{\alpha,i}$. Then $\phi(x)\Ap$ contains at the same 
time a chamber $C''$ distinct from $C$ whose closure contains $A$ and
 the half-apartment of $\Ap$ delimited by $H$ and containing $C$, which 
itself contains the closure of $C\cup\gamma(C)\cup\gamma(D)$. We deduce
 from this that we have $\gamma(\phi(x))\phi(x)\gamma(D)=\gamma(C'')$, 
hence $\phi(x)\gamma(\phi(x))D=C''$.

Moreover, $\phi(x)\gamma(\phi(x))$ is contained in $K_{C,E}$, which is a
 pro-solvable group, hence if $\gamma(\alpha)=-\alpha$, the commutator
 $[\phi(x)^{-1},\gamma(\phi(x)^{-1}]$ is an element of the subgroup $K'$
 of $K_{C,E}$ generated by $K_T$, $U_{\alpha,i+1}$ and $\gamma(U_{\alpha,i+1})$,
 which is itself contained in $K_{D\cup\gamma(D),E}$. If now 
$\gamma(\alpha)\neq\pm\alpha$ (remember that by the previous lemma we cannot
 have $\gamma(\alpha)=\alpha$), then $[\phi(x)^{-1},\gamma(\phi'x)^{-1}]$ is an
 elenent of the intersection with $K_{C,E}$ of the subgroup of $G$ generated by
 the $U_{\lambda\alpha+\mu\beta}$, where $\lambda$ and $\mu$ are positive integers
 such that $\lambda\alpha+\mu\beta$ is a root. We'll also denote by $K'$ this
 last subgroup; it is also contained in $K_{D\cup\gamma(D),E}$.

In both cases, we can apply \cite[lemma 5.1]{Cou} to see that there exists
 $k\in K'$ such that $[\phi(x)^{-1},\gamma(\phi'x)^{-1}]=\gamma(k)k^{-1}$, hence
 $\phi(x)\gamma(\phi(x))k=\gamma(\phi(x))\phi(x)\gamma(k)$. We thus have proved 
that $\phi(x)\gamma(\phi(x)k$ is an element of $K_{C,F}$ sending $D$ to $C''$; 
since this is true for any $C''$ and in particular for $C'$, $\Delta$ must contain
 all of them and $K_{C,F}$ acts transitively on them, which proves the proposition
 in this case.

Assume now that $\mathcal{C}$ does not contain $\gamma(C)$, or in other words that $C$
 and $\gamma(C)$ are separated by at least one of $H$ and $\gamma(H)$. Then they
 are separated by both of them, which means that $D$ and $\gamma(D)$ are in the 
same connected component. We can then apply the same reasoning as above with $C$
 and $D$ switched, and we obtain that for every chamber $C''$ of $X_E$ containing
 $A$ in its closure and distinct from $D$, there exists an element $g$ of $G_F$ 
such that $gC=C''$, which implies in particular that the combinatorial distance
 between $C''$ and $X_F$ must be $d$. Since by our hypothesis this is not true 
for $C'$, we must have $C'=D$ and even $\Delta=\{D\}$, and the result of the
 proposition is then trivial. $\Box$
\bigskip

\noi {\it Remark}. Actually, this very last case turns out to be impossible. 
To see that, we can for example observe that the combinatorial distance between 
$C$ and $X_F$ is equal to the combinatorial distance between $C$ and some facet 
of $X_F\cap\Ap$ of maximal dimension plus the dimension of the $F$-anisotropic 
component of $T$, and that there exists a minimal gallery between $C$ and some
 chamber of $X_F$ whose closure contains the barycenter $b$ of $C\cup\gamma(C)$ 
(which is itself an element of $X_F$); with the hypotheses of the last case, it 
is easy to check that the closure of $C\cup\{b\}$ must contain $D$, hence a contradiction.

Paul Broussous
\smallskip

paul.broussous{@}math.univ-poitiers.fr
\medskip

Fran\c cois Court\`es
\smallskip

francois.courtes{@}math.univ-poitiers.fr
\medskip

D\'epartement de Math\'ematiques

UMR 6086 du CNRS
\smallskip

T\'el\'eport 2 - BP 30179

Boulevard Marie et Pierre Curie

86962 Futuroscope Chasseneuil Cedex

 France

\end{document}